\documentclass[11pt]{article}
\title {Property P for Knots admitting Certain Gabai Disks}
\author {Oliver T. Dasbach
\thanks{e-mail: kasten@math.lsu.edu,
http://www.math.lsu.edu/$\sim$kasten}
\\Louisiana State University\\
Department of Mathematics\\Baton Rouge, LA 70803
\and
Tao Li 
\thanks{e-mail: tli@math.okstate.edu, 
http://www.math.okstate.edu/$\sim$tli}
\thanks{Partially supported by NSF DMS 0102316}
 \\ Oklahoma State University\\
Department of Mathematics\\Stillwater, OK 74078
}

\date{}

\usepackage{amsmath,amssymb,amsthm, graphics}
\usepackage{epsf,psfrag}

\newtheorem{theorem}{Theorem}[section]
\newtheorem{definition}[theorem]{Definition}

\newtheorem{lemma}[theorem]{Lemma}
\newtheorem{corollary}[theorem]{Corollary}

\newtheorem{remark}[theorem]{Remark}

\newcommand{\R} {\mathbb{R}}

\begin{document}
\maketitle
\begin{abstract}
We show that if a knot has a minimal spanning 
surface that admits certain Gabai disks, then this knot has Property P.
As one of the applications we extend and simplify a recent result of 
Menasco and Zhang that closed 3-braid knots have Property P.
Other applications are given. 
\end {abstract}

\section{Introduction}

The question whether every knot has Property P, i.e. whether non-trivial 
Dehn surgery on a knot can yield a homotopy sphere, was first raised by
Bing and Martin \cite{BingMartin:PropertyP}.
Although it was proven for many classes of knots
(e.g. \cite{DelmanRoberts:PropertyP, Scharlemann:PropertyP,
  Simon:PropertyP}), this question is still open.

Recently, Menasco and Zhang \cite{MZ:PropertyP} showed that knots represented 
by closed $3$-braids have Property P. Their interesting proof
makes use both of the Casson invariant as well as the theory of essential 
laminations. The paper of Menasco and Zhang was the main motivation for the 
results presented here.

The most striking theorem \cite{CGLS}, though, is that Property P is
``almost true'' for all knots, in the following sense: 
Let $M(r)$ be the manifold obtained by $r$-Dehn 
surgery on a knot $K$. 
Then only for $r=\pm 1$ the manifold $M(r)$ can be a 
homotopy sphere, but not both $M(-1)$ and $M(1)$ can be one.
Furthermore, by \cite {GL} it is known that non-trivial Dehn surgery 
on a knot never yields the sphere $S^3$. 
Therefore, if there is a knot without Property P then there is a 
counterexample to the Poincar\'e conjecture.

Closely related to the Property P conjecture is a slightly more general 
problem:
Is there a knot, other than the trefoil, such that Dehn surgery on it
yields a homology 3-sphere with finite fundamental group? It is well known,
that such a homology 3-sphere must be either a homotopy
3-sphere, i.e. the fundamental group is trivial, or must have as fundamental group 
the binary icosahedron group \cite{Zhang:PropertyI, Kervaire:HomologySpheres}.

The goal of this paper is to show that a non-trivial, non-trefoil knot $K$ with a minimal 
spanning surface admitting two Gabai disks, has Property P. A Gabai disk is 
an embedded disk with boundary on the spanning surface that is pierced 
transversally by the knot exactly twice and 
both times on the boundary. 
As an application we extend and re-prove the result of 
Menasco and Zhang for closed three braid knots. 

Furthermore, another immediate consequence 
is that knots represented by closed homogeneous braids have Property P. 

{\bf Acknowledgement: } Part of this work was done, while the first author was 
visiting Oklahoma State University and was enjoying the terrific hospitality 
in the Department of Mathematics. 

\section{The tools}

Two tools have been proven to be very effective in showing Property P 
for a large class of knots. The first one is the celebrated Casson invariant
for integral homology spheres (e.g. \cite{AM:Casson}). The second one
is finding essential laminations.

For an integral homology sphere the Casson invariant is defined in terms of 
``counting'' irreducible representations of the fundamental group into 
$SU(2)$. By construction, it vanishes for homotopy spheres. 

The Casson invariant has a Dehn surgery description in the following sense. 
Let $v_2(K)=1/2 \Delta''(1)(K)$,
i.e. half the second derivative of the Alexander polynomial evaluated at $1$,
be the unique (up to multiplication) Vassiliev invariant of order $2$.
The Casson invariant of the manifold $M(1/n)$ that we get by $1/n$ Dehn
surgery on $K$ is $n v_2(K)$. We will refer to $v_2(K)$ as the
Casson invariant for knots.

The Casson invariant is the first non-trivial example of 
finite type invariants for integral homology spheres (see e.g. 
\cite{Lin:FiniteTypeSurvey} for a survey on finite type invariants).
Recently, W. Li and Rubinstein \cite{LR:CassonHomotopy} proved that the Casson 
invariant is in fact a homotopy invariant. This property is not known for 
any other finite type invariant. It is even unknown, whether finite type 
invariants have to vanish on homotopy spheres.

As powerful as the Casson invariant is, it is possible
to construct an infinite family of knots such that the Casson
invariant of every homology sphere $M(1/n)$ is trivial.
Stanford gave a description of such a construction in 
\cite{Stanford}.

The second tool is based on finding essential laminations.
Informally, an essential lamination in a three manifold $M$ is a
closed subset which is foliated by suitable imbedded 2-dimensional
leaves. For a precise definition, see
\cite {GabaiOertel:EssentialLaminations}.

By a theorem of Gabai and Oertel a closed
oriented $3$-manifold $M$ which contains an essential lamination
is irreducible and has universal cover $\R^3$. Thus, in particular,
its fundamental group is not finite.

\section{Gabai Disks}

Following Ko and Lee \cite{KL:Genera} we define:

\begin{definition}
Let $S$ be a spanning surface for a knot $K$. 
An embedded disk $D$ in $S^3$ is called Gabai disk for $S$ if 
\begin{enumerate}
\item $D\cap S=\partial D$
\item $\partial D$ intersects $K$ transversally in $2$ points
\item Both arcs of $\partial D$ divided by these points are essential in $S$.
\end{enumerate}
\end{definition}

The following is a generalization of a result of Menasco and Zhang \cite{MZ:PropertyP}:

\begin{theorem} \label{main_theorem}
Let $S$ be a minimal spanning surface for the non-trivial, non-trefoil knot 
$K$.  If $(K,S)$ admits two Gabai disks that intersect in at most 
one point on $K$ then $M(r)$ contains an essential lamination for each 
rational slope $r \in (-2,2)$.
\end{theorem}

\begin{proof}
The proof is given in Section \ref{Tao}. 
\end{proof}

\section{Applications}

\subsection {Special configurations}

\begin{theorem}
Suppose, 
$K_1$ and $K_2$ are two knots that coincide outside a small ball and 
within this ball they are as in 

\medskip

\centerline{\scalebox{0.55}{\includegraphics{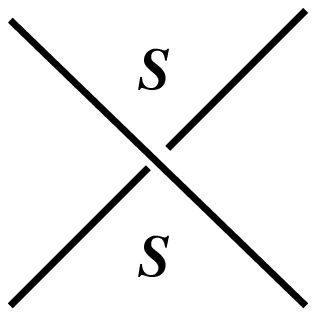}} 
\qquad \scalebox{0.55} {\includegraphics{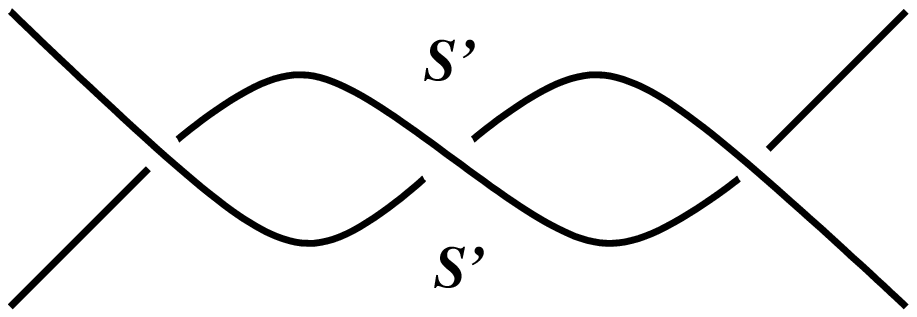}}}

and the spanning (orientable) surface $S$ of $K_1$ is minimal. 
Here, no initial assumption on the surface $S'$ are made.
Then $K_2$ has Property P.

\end{theorem}

\begin{proof}
We show that the surface $S'$ of $K_2$ is minimal, i.e. that 
$\mbox{genus}(K_2)=\mbox{genus}(K_1)+1$. 
This would follow from standard techniques using Gabai disks. 
We prefer to invoke a beautiful theorem of Scharlemann and Thompson 
\cite{ST:Linkgenus}, which was proven using the techniques of Gabai, 
though.

Suppose $L_-$ is a link obtained by changing a positive crossing in a link $L_+$ into a negative crossing and $L_0$ is the link obtained by smoothing this crossing. Then by Scharlemann's and Thompson's results two of the three numbers
$-\chi(L_+), -\chi(L_-), -\chi(L_0)+1$ are equal and the third one is not 
larger than the other two. 
Here, $-\chi(L)$ denotes the minimal negative Euler characteristic among all spanning surfaces of the link $L$.

Now let $L_{--}$ be the knot $K_2$, where the two signs stand for two of the three crossings in the diagram. Then $L_{+-}=L_{-+}=K_1=L_{00}$. 
We know that $S$ is minimal, thus $-\chi(L_{00}) > -\chi(L_{0+})$. 
This implies by the theorem of
Scharlemann and Thompson, that $$-\chi(K_1)+1=-\chi(L_{00})+1=-\chi(L_{0-}).$$
By applying the theorem for a second time, we get, as claimed, that 
$$-\chi(K_2)=-\chi(L_{--})=-\chi(L_{0-})+1=-\chi(L_{00})+2=-\chi(K_1)+2.$$

Therefore, the surface $S'$ is minimal. Since it satisfies the conditions of 
Theorem \ref{main_theorem} the claim follows.  
\end{proof}

It is interesting that a similar spirited result follows naturally by 
using the Casson invariant:

\begin{theorem}
Suppose a sequence of knots $K_n$, indexed by $n$, is defined 
by changing the number $n$ of half-twists (a negative $n$ means negative 
half-twists): 

 \medskip

\centerline{\scalebox{0.55}{\includegraphics{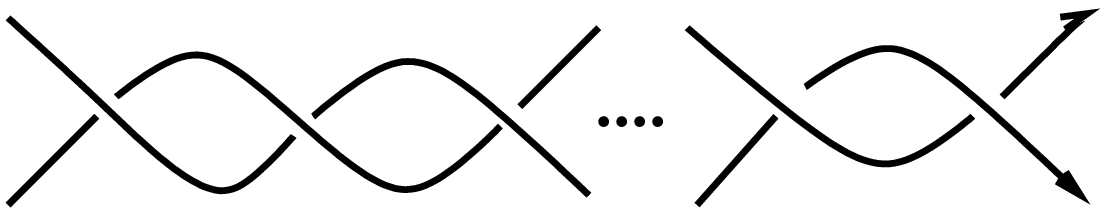}}}

Then for all but at most two values of $n$ the knot $K_n$ has Property $P$. 
\end{theorem}

\begin{proof}
If the Casson invariant $v_2(K)$ of a knot $K$ is non-zero, 
then the knot has Property $P$. Since we allow both positive and 
negative $n$ and have no restriction on the
knot outside the diagram we may assume that $n$ is odd.
For the Casson invariant it holds for $n\geq 3$):
$$v_2(T_n)-v_2(T_{n-2})=\mbox{lk}(T_{n-1})= \frac{n-1} 2 \mbox{lk} (T_0),$$
where $\mbox{lk}(L)$ denotes the linking number of a link.
A similar relation holds for negative half twists.

From this it follows for $n$ either positive or negative that 
$$v_2(T_n)= v_2(T_1) + \frac {n-1} 2 \mbox{lk}(T_0)+\frac {n^2-1} 8$$

Since this is a quadratic, non-constant polynomial in $n$ it has at most two
roots and the claim follows.
\end{proof}

\subsection {Closed $3$-braid knots}

The following theorem is a generalization and the proof is a simplification of 
a result of Menasco and Zhang in \cite{MZ:PropertyP}: 

\begin {theorem}
Let $K$ be a closed three-braid knot other than the trefoil.
Then $M(r)$ contains an essential lamination for all
$r \in (-2,2)$. In particular, $K$ has Property $P$ and no such Dehn surgery
yields a manifold with finite fundamental group.
\end{theorem}

\begin{proof} First recall the classification of closed three braids,
  given by Birman and Menasco \cite{BirmanMenasco:3braids} and the
  genus computation of these knots \cite{Xu}.
We use the presentation $$\langle a_{12},a_{23},a_{1,3} \vert a_{23}
  a_{12}= a_{13} a_{23}=a_{12} a_{13} \rangle$$ of the braid group
  $B_3$. In terms of the standard generators $\sigma_1, \sigma_2$ this means:
$a_{12}=\sigma_1, a_{23}=\sigma_2$ and 
$a_{13}=\sigma_1^{-1} \sigma_2 \sigma_1$.

Let $\beta \in B_3$ be of minimal length $l$, with respect to the generators $a_{12}, a_{13}$ and $a_{23}$, in its conjugacy class. 
It is shown in \cite{Xu} that the negative Euler characteristic of a
  minimal spanning surface $S$ of $K$ is:  $$-\chi(S)=l-3.$$ 
More specifically, the following generalized Seifert 
algorithm produces a minimal spanning surface: 
The surface is constructed by three disks and each 
$a_{ij}$ represents a band that connects disk $i$ with disk $j$.

We will show that every three strand braid $\beta$ which represents a
knot $K$ either contains, up to length preserving conjugacy and braid isotopy,
two Gabai disks as in Theorem \ref{main_theorem} 
or $K$ is the trefoil or the knot $5_2$. 

From the presentation for $B_3$ we get the following relations:
\begin{eqnarray*}
a_{12} a_{23}^{-1}=a_{23}^{-1} a_{13}, &\qquad& a_{12}
a_{13}^{-1}=a_{23}^{-1} a_{12}\\
a_{23} a_{12}^{-1}= a_{13}^{-1} a_{23} &\qquad& a_{23} a_{13}^{-1} =
a_{13}^{-1} a_{12} \\
a_{13} a_{12}^{-1} = a_{12}^{-1} a_{23} &\qquad& a_{13} a_{23}^{-1} =
a_{12}^{-1} a_{13}.
\end{eqnarray*}

Thus, we can write every braid $\beta \in B_3$ as $P N$ for some
positive word $P$ and negative word $N$ without increasing its length.

If $PN$ has minimal length in its conjugacy class, we can find Gabai disks 
as in Theorem \ref{main_theorem} in the following way: 
If the exponent sum in P of one of the three generators is greater than 
1 than we have one Gabai disk. 
If it is greater than two, we have two such Gabai disks. 
Accordingly if the exponent sum in N is smaller than $-1$ than we have one 
Gabai disk, if it is at most $-3$ we have two Gabai disks.

This already reduces our considerations to a finite number of knots 
that we have to check: The length of both $P$ and $N$ has to be $\leq 4$. 
The rest of the proof is concerned with this finite number. 
It could be easily verified on a computer. 

With $\delta=a_{23} a_{12}=a_{13} a_{23}=a_{12} a_{13}$
we have
\begin{equation} \label{delta}
\delta a_{23}=a_{12} \delta, \qquad \delta a_{13} = a_{23}
\delta, \qquad \delta a_{12}=a_{13} \delta.
\end{equation}
Therefore we can assume that $\beta = \delta^{k} P N$ for some $k$ and
$P$ and $N$ do not contain $\delta$ or $\delta^{-1}$ as a subword.

Since $\delta$ starts and ends with every generator and $\beta$ is assumed to
be minimal in its conjugacy class we can assume that either $k = 0$ 
or $k<0$ and $P$ is empty or $k>0$ and $N$ is empty.

If $k\geq 2$ or $k \leq -2$ we have two Gabai disks. 

Thus, to summarize, $\beta$ is either $\delta P$, $P N$ or $\delta^{-1} N$
where $P$ and $N$ do not contain powers of $\delta$ as a subword and the
lengths of both $P$ and $N$ is $\leq 4 $. 
Moreover, since we assume a knot, the length of $\beta$ is even and since 
the knot is non-trivial, its length must be $\geq 4$.

Up to conjugation we can assume that a non-empty $P$ is one of the following 
seven words:
$a_{12}, a_{12}^2, a_{12} a_{23}, a_{12}^2 a_{23}, a_{12} a_{23} a_{13}, 
a_{12}^2 a_{23} a_{13}, a_{12} a_{23} a_{13} a_{12}$.

In particular, every $\delta P$ of length at least 4 admits at least two 
Gabai disks.
A corresponding argument shows that every $\delta^{-1} N$ admits at
least two Gabai disks.

Finally, assume that $\beta = P N$ and $\beta$ does not contain two 
Gabai disks.
Conjugation with $\delta$ (see Equation (\ref{delta})) induces an inner
automorphism in $B_3$ that sends a generator $a_{ij}$ to any other of the
other two generators. In particular, we can assume that $P$ starts 
with $a_{12}$ and, thus, $N$ does not end with $a_{12}^{-1}$.

Assume first, that the length of $P$ is greater or equal the length of $N$.
Moreover, if $P$ and $N$ are of equal length, then $N$ does not contain
the square of a generator as a subword. 
Since not both $P$ and $N$ can have length $4$, we can assume that the
length of $N$ is at most $3$. 
If $N$ is empty, then $P$ must be $a_{12}^2 a_{23} a_{13}$, i.e. the knot 
$5_2$.  

If the length of $N$ is one, then the length of $P$ must be $3$ and 
we have the possibilities:
$a_{12}^2 a_{23} a_{13}^{-1}=a_{12}^2 a_{13}^{-1} a_{12}$ and we have 
two Gabai disks. Otherwise, it is $a_{12} a_{23} a_{13} a_{23}^{-1}$ 
which is a $3$-component link.

If the length of $N$ is two then the length of $P$ is either $2$ or $4$.
The cases are: $a_{12}^2 a_{13}^{-1} a_{23}^{-1}=a_{12} a_{23}^{-1} a_{12} a_{23}^{-1}, a_{12}^2 a_{23} a_{13} a_{23}^{-1} a_{12}^{-1}=a_{12}^2 a_{23} a_{12}^{-1} a_{13} a_{12}^{-1}, a_{12}^2 a_{23} a_{13} a_{12}^{-1} a_{13}^{-1}=a_{12}^2 a_{23} a_{12}^{-1} a_{23} a_{13}^{-1}$ and $a_{12} a_{23} a_{13} a_{12} a_{13}^{-1} a_{23}^{-1}=a_{12} a_{23} a_{13} a_{23}^{-1} a_{12} a_{23}^{-1}$.
All of these knots contain two Gabai disks. 
Furthermore, $a_{12} a_{23} a_{12}^{-1} a_{13}^{-1}$ 
and $a_{12} a_{23} a_{13}^{-1} a_{23}^{-1}$ that are three component links.

If the length of both $P$ and $N$ is three then we have the three
possibilities, all of them containing two Gabai disks:
$a_{12}^2 a_{23} a_{12}^{-1} a_{13}^{-1} a_{23}^{-1} = a_{12}^2 a_{13}^{-1} a_{23} a_{13}^{-1} a_{23}$, $a_{12} a_{23} a_{13} a_{12}^{-1} a_{13}^{-1} a_{23}^{-1} = a_{12} a_{13}^{-1} a_{23}^2 a_{13}^{-1} a_{23}^{-1}$ and $a_{12} a_{23} a_{13} a_{23}^{-1} a_{12}^{-1} a_{13}^{-1} = a_{12} a_{23} a_{12}^{-2} a_{23} a_{13}^{-1}$.

Finally, if the length of $N$ is greater or equal the length of $P$ then a
case by case check gives us the only exception for a knot without two Gabai 
$2$-disks to be $a_{12}^{-1} a_{13}^{-1} a_{23}^{-1}$ which is again the knot
$5_2$.

The final case, i.e. the knot $5_2$, is covered by a result of Delman and Roberts
\cite{DelmanRoberts:PropertyP}. Suppose $K$ is a non-torus alternating knot
such that the checkerboard surface is the one that one gets by applying the 
Seifert algorithm to its knot diagram. Then every finite Dehn surgery on 
$K$ produces a manifold containing an essential lamination.
\end{proof}

\subsection {Knots represented by closed homogeneous braids}

One class of knots seems to be especially made for the application of
Theorem \ref{main_theorem}: Knots represented by closed homogeneous 
braids. A homogeneous braid is a braid where every generator $\sigma_i$
either occurs always with positive or always with negative exponent.

\begin{theorem}
Property P holds for non-trivial knots represented by 
closed homogeneous braids.
Furthermore, unless the knot is the trefoil, no Dehn surgery along 
the knot yields a homology sphere with finite fundamental group.
\end{theorem}

\begin{proof}
By a result of Stallings \cite{Stallings:fibred_links} a knot represented by a 
closed homogeneous braid is fibred and a minimal spanning surface $S$ 
is given by Seifert's algorithm. Since the braid is homogeneous and 
non-trivial, the surface $S$ admits at least two Gabai disks.
\end{proof}

\section{Essential laminations and Gabai disks}
\label{Tao}

Essential laminations, which are introduced by Gabai and Oertel \cite{GabaiOertel:EssentialLaminations}, are a generalization of both incompressible surfaces and taut foliations.  We first give a brief review on the basic definitions and properties of essential laminations and we refer the reader to \cite{GabaiOertel:EssentialLaminations, TaoLi:Laminar, TaoLi:TrainTracks} for details.

\begin{definition}[\cite{GabaiOertel:EssentialLaminations}]
Let $N$ be a close 3-manifold and $\lambda$ be a lamination in $N$.  We say $\lambda$ is an \emph{essential lamination} in $N$ if it satisfies the following conditions. 
\begin{enumerate}
    \item The inclusion of leaves of $\lambda$ into $N$ induces an injection on $\pi_1$.
    \item The complement of $\lambda$ is irreducible.
    \item $\lambda$ has no sphere leaves.
    \item $\lambda$ is end-incompressible. 
\end{enumerate}
\end{definition}

A closed 3-manifold is said to be \emph{laminar} if it contains an essential lamination. It is shown in \cite{GabaiOertel:EssentialLaminations} that the universal cover of a laminar 3-manifold is $\mathbb{R}^3$.  Moreover, it follows from \cite{GK:Laminations, Calegari:RFoliations, Calegari:Promoting} that laminar 3-manifolds satisfy the weak hyperbolization conjecture, i.e. either the fundamental group contains a $\mathbb{Z}\oplus\mathbb{Z}$ subgroup, or the manifold has word-hyperbolic fundamental group.  

Branched surfaces, which are a generalization of train tracks, are very useful in the theory of essential laminations and incompressible surfaces, see \cite{TaoLi:Laminar} for details.  

Let $B$ be a branched surface and $L$ be the branch locus of $B$.  We call the closure (under path metric) of each component of $B-L$ a \emph{sector} of $B$.  $L$ is a collection of smooth immersed curves in $B$.  Let $Z$ be the union of double points of $L$.  We associate with every component of $L-Z$ a normal vector (in $B$) pointing in the direction of the cusp.  We call it the \emph{branch direction} of this arc \cite{TaoLi:Laminar}.  We call a disk sector of $B$ a \emph{sink disk} if the branch direction of every smooth arc (or curve) in its boundary points into the disk.  Let $N(B)$ be a regular neighborhood of $B$.  $N(B)$ can be considered as an $I$-bundle over $B$. The boundary of $N(B)$ consists of $\partial_hN(B)$ and $\partial_vN(B)$, which are called the horizontal boundary and the vertical boundary respectively.  $B$ is said to be incompressible if 
\begin{enumerate}
\item $\partial_hN(B)$ is incompressible in $N-int(N(B))$, $\partial_hN(B)$ has no sphere component, and $N-int(N(B))$ is irreducible; 
\item $N-int(N(B))$ contains no monogon; 
\item $N(B)$ contains no disk of contact.  
\end{enumerate}
Suppose $B$ is incompressible.  A \emph{trivial bubble} is a 3-ball in $N-B$ that can be trivially eliminated by pinching $B$ (see \cite{TaoLi:Laminar} for more details).

\begin{definition}[\cite{TaoLi:Laminar, TaoLi:TrainTracks}]\label{D:lbs}
Let $B$ is a properly embedded branched surface in a 3-manifold $M$.  We say $B$ is a \emph{laminar branched surface} if, after eliminating trivial bubbles, the followings are satisfied.
\begin{enumerate}
\item$\partial_hN(B)$ is essential in the following sense: $\partial_hN(B)$ is incompressible and $\partial$-incompressible in $M-int(N(B))$, there is no monogon in $M-int(N(B))$ and no component of $\partial_hN(B)$ is a sphere or a disk properly embedded in $M$.
\item $M-int(N(B))$ is irreducible and $\partial M-int(N(B))$ is incompressible in $M-int(N(B))$.
\item $B$ contains no Reeb branched surface (see \cite{GabaiOertel:EssentialLaminations} for the definition of Reeb branched surface).
\item $B$ has no sink disk or half sink disk.
\end{enumerate}
\end{definition}

The following theorem gives a certain equivalence relation between essential laminations and branched surfaces without sink disks. 

\begin{theorem}[\cite{TaoLi:Laminar}]\label{TL}
Suppose $N$ is a closed and orientable 3-manifold.  Then
\begin{enumerate}
\item [(a)] Every laminar branched surface in $N$ fully carries an essential lamination.

\item [(b)] Any essential lamination in $N$ that is not a lamination by planes is fully carried by a laminar branched surface. 
\end{enumerate}
Furthermore, if $\lambda\subset N$ is a lamination by planes (hence $N=T^3$), then any branched surface carrying $\lambda$ must contain a sink disk and hence is not a laminar branched surface.
\end{theorem}

\begin{figure}[h]
\begin{center}
\resizebox{4in}{!}{\includegraphics{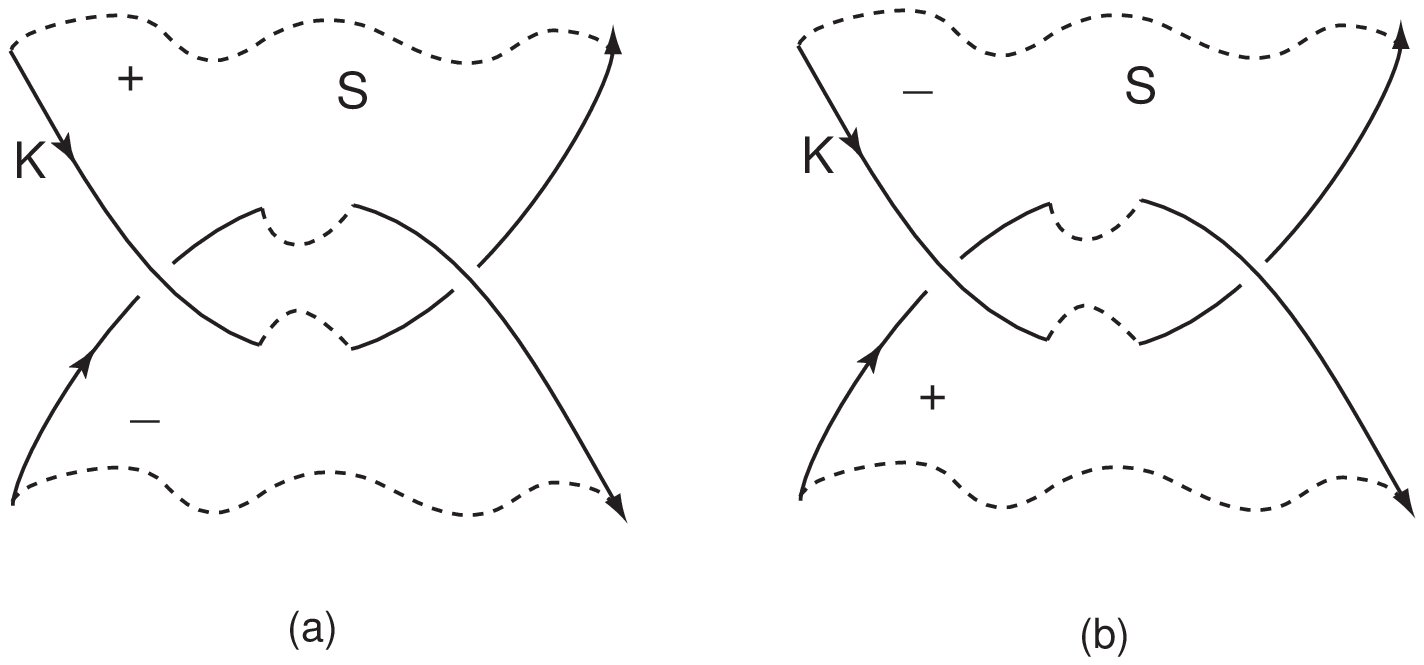}}
\caption{}\label{twists}
\end{center}
\end{figure}

A relative version of Theorem~\ref{TL} for knot manifolds is proved in \cite{TaoLi:TrainTracks}.
Let $M$ be an irreducible 3-manifold whose boundary is an incompressible torus.  Suppose $B$ is a properly embedded branched surface with boundary ($\partial B\subset\partial M$).  Let $L$ be the branch locus of $B$, and $D$ be the closure (in the path metric) of a disk component of $B-L$.  We call $D$ a \emph{half sink disk} if $\partial D\cap\partial M\ne\emptyset$ and the branch direction of each arc in $\partial D-\partial M$ points into $D$.  Note that $\partial D\cap\partial M$ may not be connected.

Let $T$ be a torus, then the isotopy class of every nontrivial simple closed curve in $T$ uniquely corresponds to a rational slope in $\mathbb{Q}\cup\infty$.  Let $\tau$ be a train track in $T$ and $s\in\mathbb{Q}\cup\infty$ be a rational slope.  We say that $s$ can be \emph{fully carried} by $\tau$ (or \emph{realized} by $\tau$) if $\tau$ fully carries a union of simple closed curves with slope $s$, i.e., we can split $\tau$ into  a union of simple closed curves with slope $s$.  For any slope $s\in\mathbb{Q}\cup\infty$, we denote by $M(s)$ the manifold after Dehn filling along $s$.

\begin{theorem}[\cite{TaoLi:TrainTracks}]\label{T:boundary}
Let $M$ be an irreducible and orientable 3-manifold whose boundary is an incompressible torus.  Suppose $B$ is a laminar branched surface and $\partial M-\partial B$ is a union of bigons.  Then, for any rational slope $s\in\mathbb{Q}\cup\infty$ that can be fully carried by the train track $\partial B$, if $B$ does not carry a torus that bounds a solid torus in $M(s)$, $M(s)$ contains an essential lamination.
\end{theorem}
\begin{remark}
Note that since $B$ is incompressible, $\partial M-\partial B$ is a union of bigons and annuli.  So, if $\partial B$ fully carries more than one slope, $\partial M-\partial B$ is always a union of bigons.  Moreover, if $B$ does not carry any closed surface, then $B$ does not carry any torus and the last requirement in the theorem is satisfied.
\end{remark}

In this section, we will mainly apply Theorem~\ref{T:boundary} to certain knot complements.  Let $K$ be a knot in $S^3$, and $M=S^3-int(N(K))$, where $N(K)$ is a tubular neighborhood of $K$.  Suppose $S$ is a minimal genus Seifert surface of $K$.  To simplify notation, we also use $S$ to denote the corresponding surface properly embedded in $M$ (i.e. $S-int(N(K)$).  We first fix an orientation for $K$ and assign a ``$+$" and ``$-$" side to $S$.

\begin{figure}[h]
\begin{center}
\resizebox{4in}{!}{\includegraphics{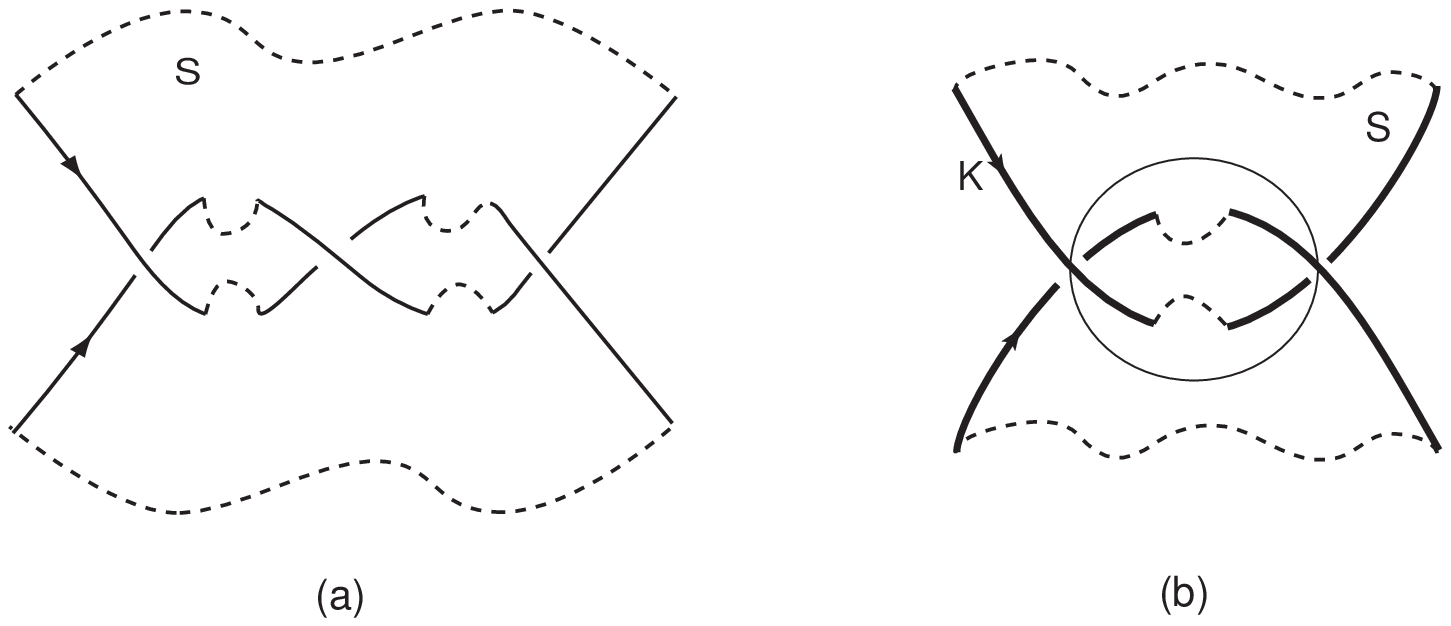}}
\caption{}\label{triple}
\end{center}
\end{figure}

\begin{definition}\label{newdef}
We say $S$ has a positive (resp. negative) \emph{double twist}, if there is a 3-ball whose intersection with $K$ and $S$ is as shown in Figure~\ref{twists}(a) (resp. (b)), and inside the 3-ball, one can add a Gabai disk (above Figure~\ref{twists}) as shown in Figure~\ref{triple}(b).  Note that the dashed curves in Figure~\ref{twists} (and other figures in this paper) denote the intersection of $S$ and the boundary of the 3-ball.  The circle in Figure~\ref{triple}(b) denote the circle along which we add the Gabai disk.  Similarly, we say $S$ has a (positive or negative) \emph{triple twist} if there is a 3-ball whose intersection with $K$ and $S$ is as shown in Figure~\ref{triple}(a), and one can simultaneously add two Gabai disks (both above Figure~\ref{triple}(a)) in the 3-ball.  Note that any two crossings in Figure~\ref{triple}(a) form a (positive or negative) double twist, and one can add a Gabai disk (from above) along any two crossings.  
\end{definition}

If $S$ has a double twist and one adds a Gabai disk as in Figure~\ref{triple}(b) (see \cite{Gabai:GeneraAlternating, Gabai:Fibred}), by fixing a normal direction for the Gabai disk, one can obtain a taut sutured manifold decomposition for the sutured manifold $\overline{M-S}$ \cite{Gabai:Fibred, Gabai:FoliationsI}.  Note that the Gabai disk corresponds to a product disk in $\overline{M-S}$ and one always gets a taut sutured manifold by adding any product disk to a taut suture manifold \cite{Gabai:Fibred, Gabai:FoliationsI}.  Let $D$ be the product disk in $\overline{M-S}$ corresponding to the Gabai disk.  Then, one can modify $S\cup D$ into a branched surface $B$ properly embedded in $M$, and $M-int(N(B))$ is the taut sutured manifold above.  The branched locus of $B$ consists of two disjoint arcs properly embedded in $S$.  Figure~\ref{traintrack1}(a) is a schematic picture of the branch locus of $B$ (the arrows denote the branch directions).  The boundary train track $\partial B$ is as shown in Figure~\ref{traintrack1}(b), see Figure 14 of \cite{RRoberts:Constructing} for a picture of $B$ near $\partial M$.  

\begin{figure}[h]
\begin{center}
\resizebox{4in}{!}{\includegraphics{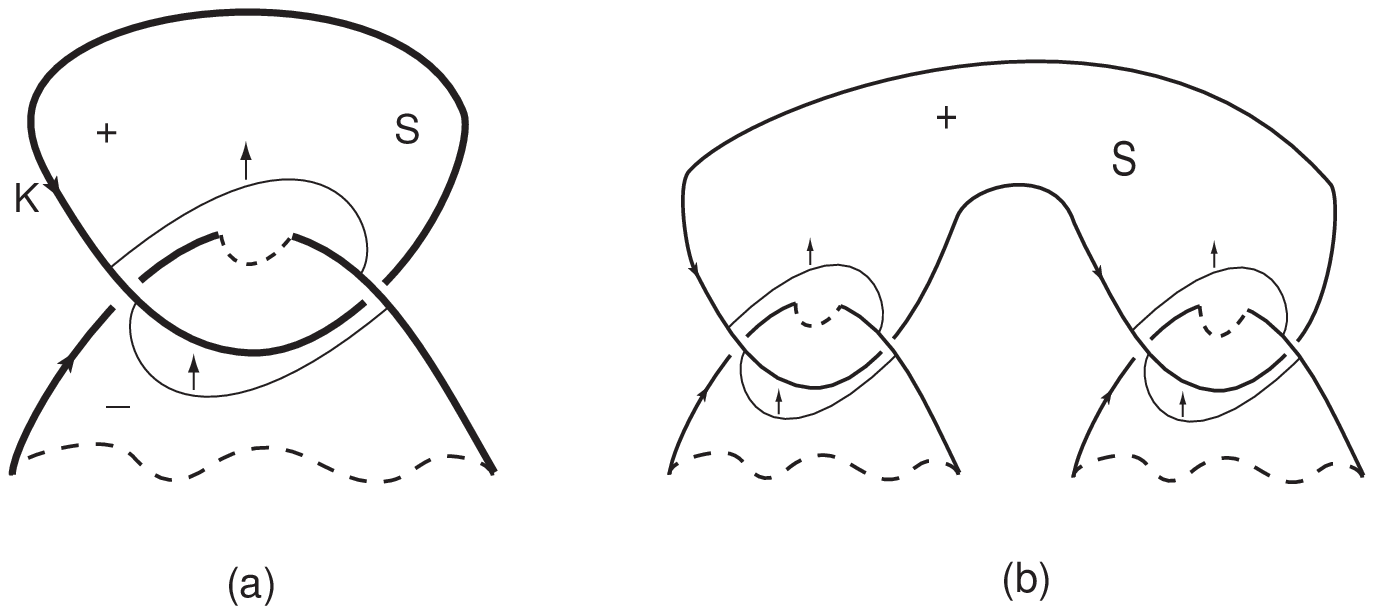}}
\caption{}\label{disk}
\end{center}
\end{figure}

\begin{figure}[h]
\begin{center}
\resizebox{3in}{!}{\includegraphics{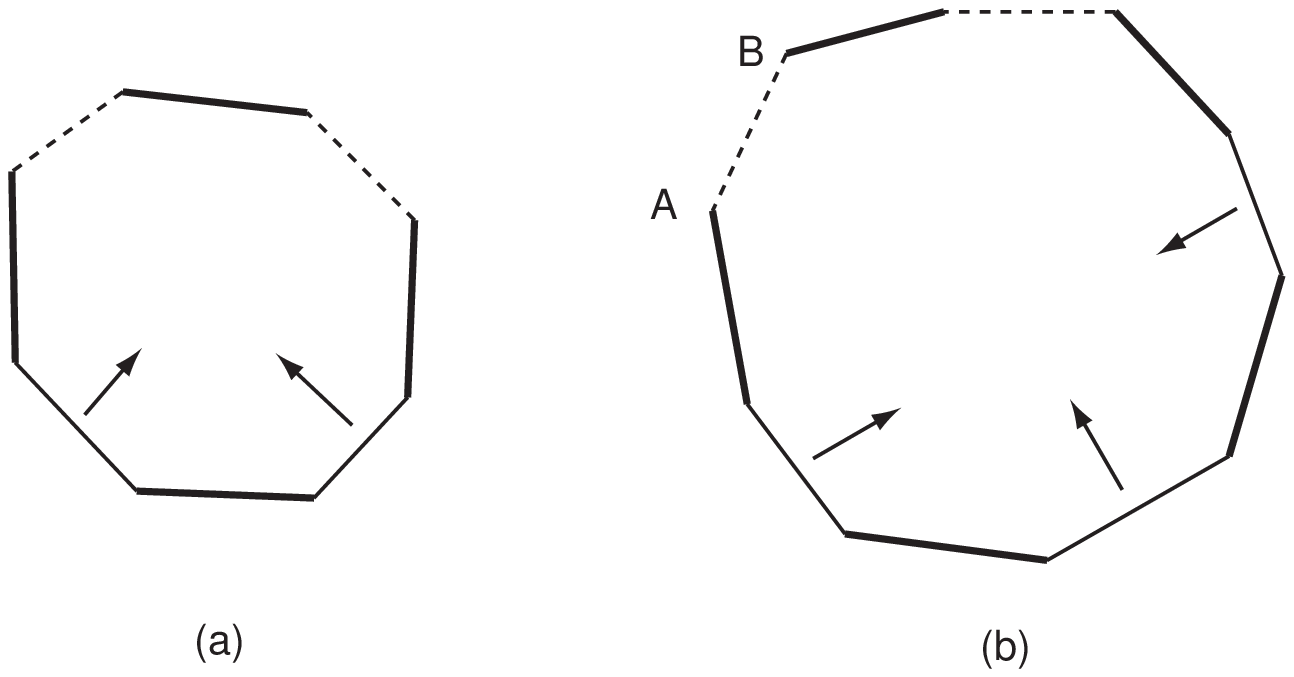}}
\caption{}\label{sink}
\end{center}
\end{figure}

\begin{figure}[h]
\begin{center}
\resizebox{3.5in}{!}{\includegraphics{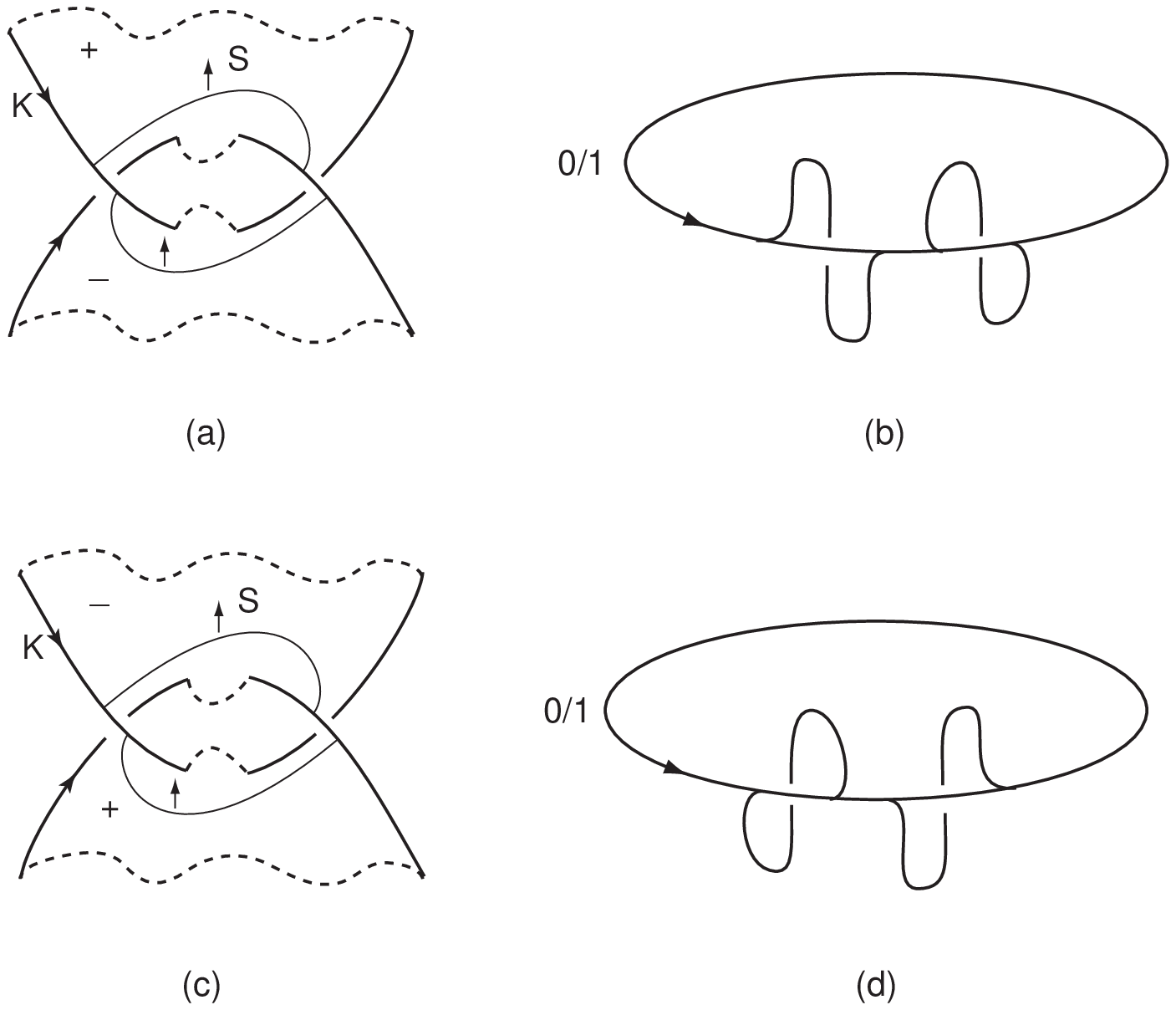}}
\caption{}\label{traintrack1}
\end{center}
\end{figure}

\begin{lemma}[\cite{RRoberts:TautFoliation} and compare with \cite{MZ:PropertyP}]
\label{twist1}
If $S$ has a positive (resp. negative) double twist, then $M(r)$ contains a taut foliation for any $r\in (-\infty,1)$ (resp. $r\in (-1, \infty)$). 
\end{lemma}
\begin{proof}
The proof is similar to \cite{RRoberts:TautFoliation} in spirit.  Suppose $S$ has a positive double twist.  Let $B$ be the branched surface above.  We first show that $B$ contains no (half) sink disk.  The branch sector corresponding to the Gabai disk $D$ is not a half sink disk as the branch directions all point out of $D$.  Now, we consider the local picture of $S$ in Figure~\ref{traintrack1}(a).  The subsurface of $S$ in Figure~\ref{traintrack1}(a) is an annulus and the branch locus consists of two essential arcs in this annulus that cut this annulus into two disks, say $d_1$ and $d_2$.  The branch directions all point into $d_1$ and out of $d_2$.  Figure~\ref{sink}(a) is a picture of $d_1$ (the arcs with arrows denote the branch locus, the dashed arcs denote the intersection of $d_1$ with the boundary of the 3-ball in definition~\ref{newdef}, and the other arcs are from $\partial S$).  Since any branch sector that does not contain $d_1$ has a boundary arc with branch direction pointing outwards, if there is a (half) sink disk, then the two dashed lines in Figure~\ref{sink}(a) must be parallel to $K$ (in $S$), and hence the picture must be as shown in Figure~\ref{disk}(a), in which case $K$ is a link with more than one component.  This contradicts our hypothesis that $K$ is a knot.  Therefore, $B$ cannot have any (half) sink disk.

Since $M-int(N(B))$ is a taut sutured manifold and the suture is $\partial_vN(B)\cup (\partial M-int(N(B)))$ \cite{Gabai:FoliationsI, Gabai:FoliationsIII}, it follows from our construction that $B$ is a laminar branched surface.  Moreover, $\partial M-\partial B$ is a union of bigons and $B$ does not carry any closed surface.  Therefore, $B$ satisfies all the  hypotheses in Theorem~\ref{T:boundary}, and hence $M(r)$ contains an essential lamination for every slope fully carried by $\partial B$.  If $S$ has a positive double twist, then the train track $\partial B$ (corresponding to Figure~\ref{traintrack1}(a)) is as shown in Figure~\ref{traintrack1}(b).  By \cite{TaoLi:TrainTracks, RRoberts:TautFoliation} (or a simple calculation), the interval of slopes fully carried by this train track is $(-\infty, 1)$

If $S$ has a negative double twist, then the boundary train track $\partial B$ (corresponding to Figure~\ref{traintrack1}(c)) is as shown in Figure~\ref{traintrack1}(d), and  the interval of slopes fully carried by this train track is $(-1,\infty)$.

Moreover, since $B$ is obtained by a sutured manifold decomposition and $K$ is a knot in $S^3$, by \cite{Gabai:FoliationsIII}, the essential laminations above can be extended to taut foliations in $M(r)$.
\end{proof}

\begin{corollary}\label{both}
If $S$ contains both positive and negative double twists, then $M(r)$ contains a taut foliation for any slope $r\ne\infty$.
\end{corollary}
\qed

\begin{figure}[h]
\begin{center}
\resizebox{4.5in}{!}{\includegraphics{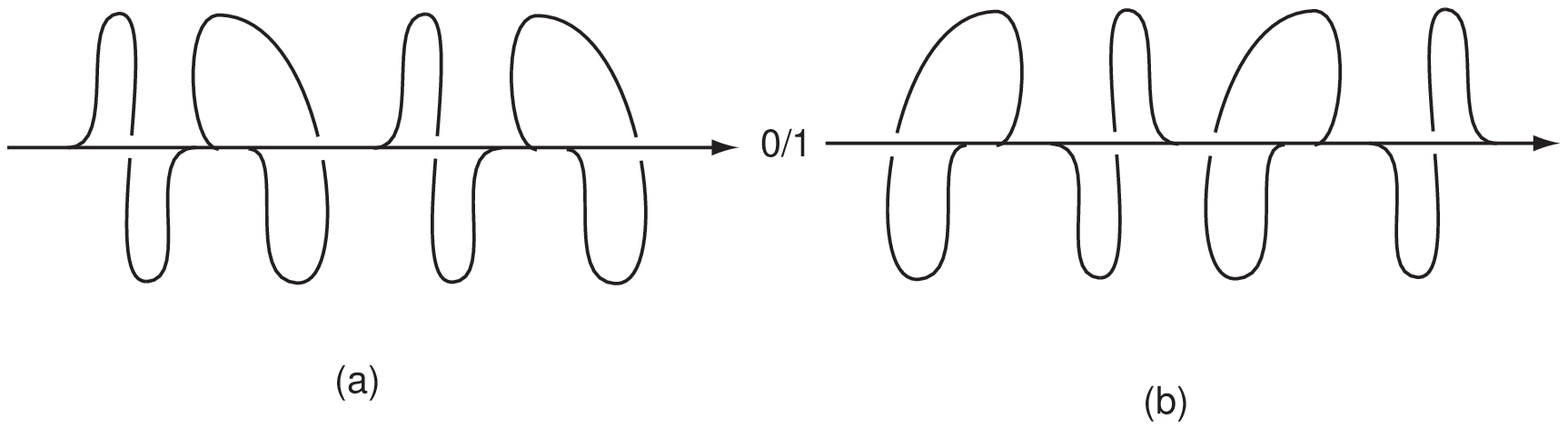}}
\caption{}\label{traintrack2}
\end{center}
\end{figure}

\begin{lemma}\label{twist2}
If $S$ has two disjoint positive (resp. negative) double twists, in other words, if there are two disjoint 3-balls whose intersections with $S$ and $K$ are as shown in Figure~\ref{twists}(a) (resp. (b)), then $M(r)$ contains a taut foliation for any slope $r\in (-\infty, 2)$ (resp. $(-2,\infty)$).
\end{lemma}
\begin{proof}
If there are two disjoint positive double twists, then we can add two Gabai disks $D_1$ and $D_2$, which give rise to a taut sutured manifold decomposition. As before, we deform $S\cup D_1\cup D_2$ into a branched surface $B$ so that $M-int(N(B))$ is a taut sutured manifold.  Since we have two Gabai disks, the boundary train track $\partial B$ is as shown in Figure~\ref{traintrack2}(a), and it follows from \cite{TaoLi:TrainTracks} that interval of slopes fully carried by such a train track is $(-\infty, 2)$.  If $S$ has two disjoint negative twists, then the corresponding boundary train track is as shown in Figure~\ref{traintrack2}(b), in which case the interval of slopes in $(-2,\infty)$.

Similar to the argument in the proof of Lemma~\ref{twist1}, we have two disjoint 3-balls, which together with the branch locus cut $S$ into pieces.  In this case, we have two copies of Figure~\ref{sink}(a), and any other piece contains a boundary arc with branch direction pointing outwards.  If $B$ contains a (half) sink disk, then either $K$ and $S$ have a local picture as in Figure~\ref{disk}(a) and Lemma~\ref{twist1}, or two dashed arcs from the two copies of Figure~\ref{sink}(a) are parallel in $S$ and the other two dashed arcs are parallel to $K$ in $S$.  Note that when we connect the two copies of Figure~\ref{sink}(a) along dashed arcs, we have to make the $+$ and $-$ sides compatible.  Figure~\ref{disk}(b) gives an example of getting a (half) sink disk by connecting the two copies of Figure~\ref{sink}(a), and there are other possible configurations.

However, it is easy to see that, in any cases, if one obtains a (half) sink disk by connecting the dashed arcs in the two copies of Figure~\ref{sink}(a), $K$ is always a link with more than one component, as illustrated in Figure~\ref{disk}(a) and (b).  Thus, as in Lemma~\ref{twist1}, $B$ is a laminar branched surface and the lemma follows from Theorem~\ref{T:boundary}.  
\end{proof}

\begin{figure}[h]
\begin{center}
\resizebox{4in}{!}{\includegraphics{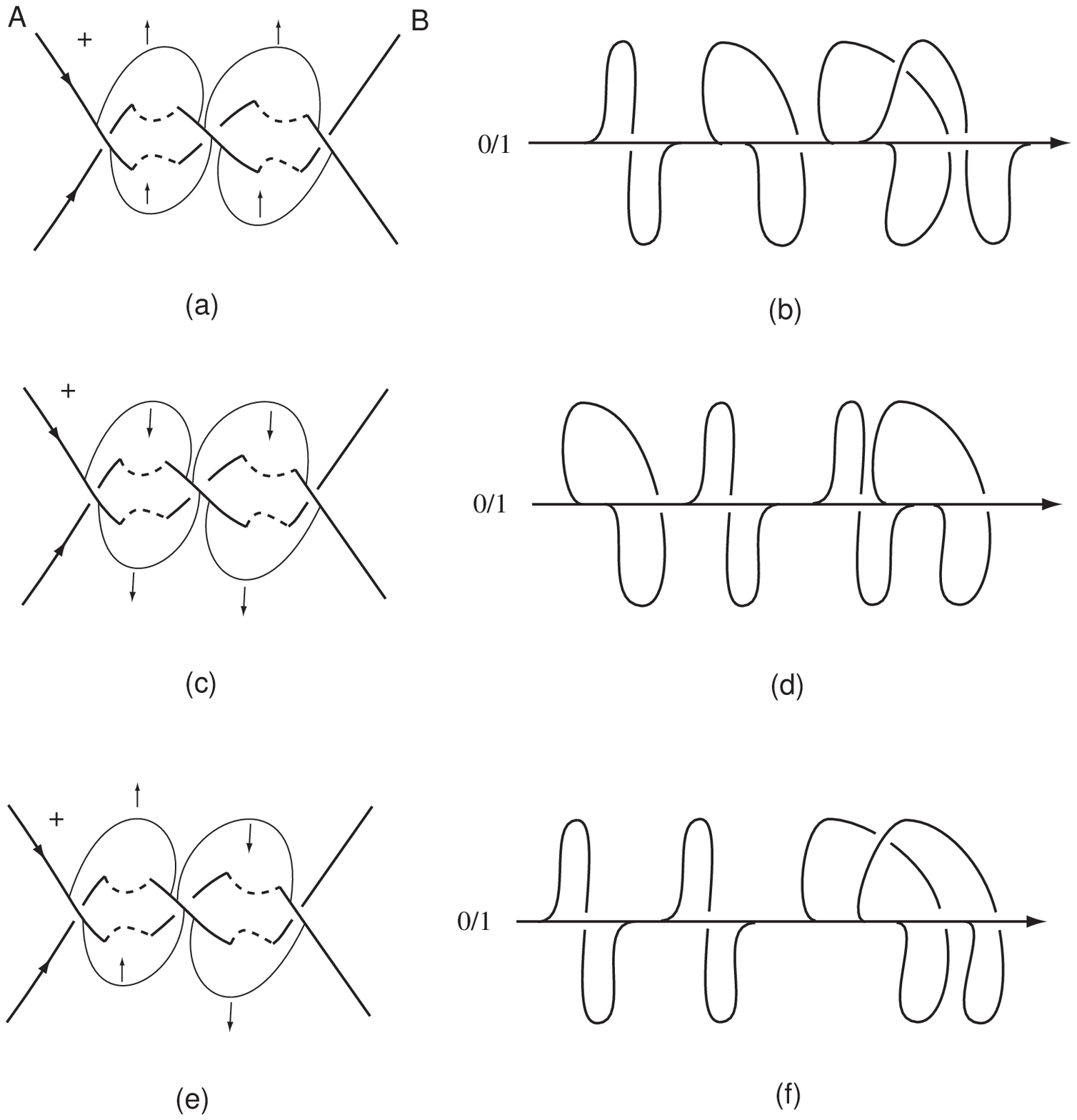}}
\caption{}\label{traintrack3}
\end{center}
\end{figure}

\begin{figure}[h]
\begin{center}
\resizebox{3.7in}{!}{\includegraphics{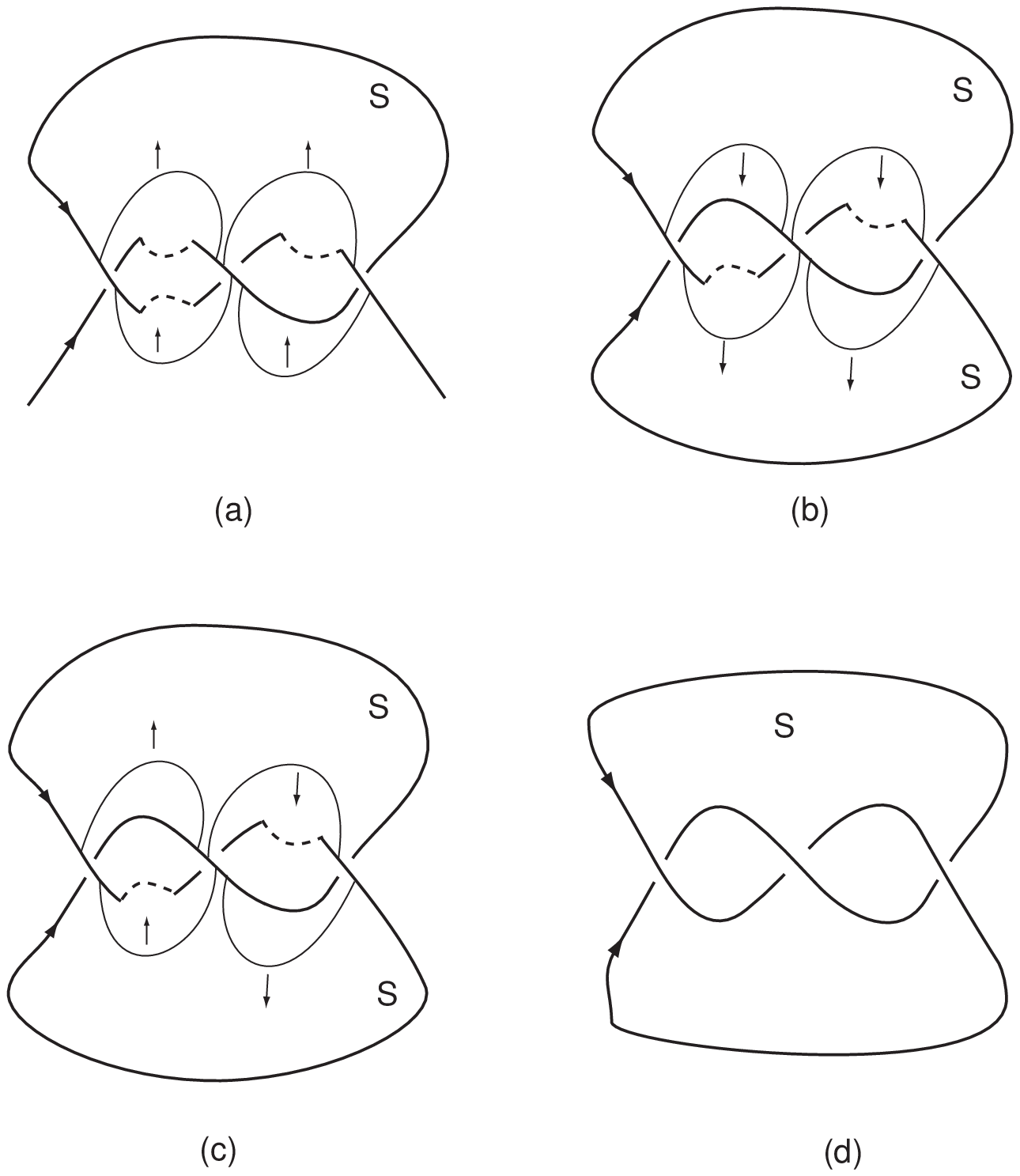}}
\caption{}\label{bad}
\end{center}
\end{figure}

\begin{lemma}\label{twist3}
Suppose $K$ is not a trefoil knot and $S$ has a positive (resp. negative) triple twist.  Then, $M(r)$ contains a taut foliation for any $r\in (-\infty, 2)$ (resp. $(-2,\infty)$).
\end{lemma}
\begin{proof}
We will only prove the case that the triple twist is positive, and the negative case is similar.  If we have a triple twist, we can still add two Gabai disks, i.e., add two product disks to $\overline{M-S}$, and get a branched surface $B$ such that $M-int(N(B))$ is the corresponding taut sutured manifold.  The schematic picture of the branch locus of $B$ in this case is as shown in Figure~\ref{traintrack3}(a), and the corresponding boundary train track $\partial B$ is as shown in Figure~\ref{traintrack3}(b).  

If we reverse the orientation of both Gabai disks, we will get a different branched surface $B'$ but $M-int(N(B'))$ is the same taut sutured manifold as $M-int(N(B))$ \cite{Gabai:Fibred, Gabai:FoliationsI, Gabai:FoliationsIII}.  Figure~\ref{traintrack3}(c) is a schematic picture of the branch locus of $B'$, and the boundary train track $\partial B'$ is as shown in Figure~\ref{traintrack3}(d).  Furthermore, if we reverse the orientation of only one Gabai disk, as shown in Figure~\ref{traintrack3}(e), then we get a branched surface $B''$ which also gives rise to the same taut sutured manifold decomposition \cite{Gabai:Fibred, Gabai:FoliationsI, Gabai:FoliationsIII}.  The boundary train track $\partial B''$ is as shown in Figure~\ref{traintrack3}(f).  It following from \cite{TaoLi:TrainTracks} that the interval of realizable slopes for the three train tracks $\partial B$, $\partial B'$ and $\partial B''$ are the same, i.e., $(-\infty, 2)$.

Similar to the proof of Lemma~\ref{twist1}, we consider the local picture of $S$, i.e. the intersection of $S$ with the 3-ball in definition~\ref{newdef}.  The branch locus of $B$ cuts the subsurface in Figure~\ref{traintrack3}(a) into pieces and only one piece, say $d$, has all the branch directions pointing inward.  Figure~\ref{sink}(b) is a picture of $d$, where the arcs with arrows denote the branch locus, the dashed lines denote the intersection of $d$ with the boundary of the 3-ball, all other boundary arcs are from $K$, and the two points $A$, $B$ denote the points $A$, $B$ in Figure~\ref{traintrack3}(a).  Similar to the proof of Lemma~\ref{twist1}, if $B$ contains a (half) sink disk, then the two dashed lines in Figure~\ref{sink}(b) must be parallel to $K$ in $S$, and hence the configuration of $K$ must be as shown in Figure~\ref{bad}(a).

We suppose $B$ contains a sink disk and hence the configuration of $K$ is as in Figure~\ref{bad}(a).  Then, we reverse the orientation of both Gabai disks and construct another branched surface $B'$ as above.  Similarly, if $B'$ also contains a sink disk, then the configuration of $K$ must be as in Figure~\ref{bad}(b).  Now, we suppose both $B$ and $B'$ contain sink disks, and hence the configuration of $K$ is as in Figure~\ref{bad}(b).  Then we reverse the orientation of only one Gabai disk, as in Figures \ref{traintrack3}(e) and \ref{bad}(c), and construct the branched surface $B''$ as above. By the same argument, if $B''$ also contains a sink disk, then the two dashed lines in Figure~\ref{bad}(c) must be parallel to $K$ in $S$, and hence the configuration of $K$ must be Figure~\ref{bad}(d) which is a trefoil knot.

Since we have assumed that $K$ is not a trefoil knot, at least one of the three branched surfaces $B$, $B'$ and $B''$ does not contain any sink disk, and the lemma follows from Theorem~\ref{T:boundary} and the discussion of the train tracks of the three branched surfaces above.
\end{proof}

The next theorem is an immediate corollary of Lemmas~\ref{twist1}, \ref{twist2}, \ref{twist3} and Corollary~\ref{both}.

\begin{theorem} \label{Theorem511}
Let $K$ be a knot in $S^3$ and $S$ be a minimal genus Seifert surface.  Suppose $S$ has a triple twist or two disjoint double twists.  Then, unless $K$ is a trefoil knot, any Dehn surgery along the slope $r\in(-2,2)$ yields a 3-manifold with infinite fundamental group.  In particular, $K$ has property P. Moreover, if a Dehn surgery yields a Poincare homology sphere, then $K$ must be a trefoil knot and the slope is $1$ or $-1$.
\end{theorem}
\begin{proof}
By a simple homology argument, the Dehn surgery along the slope $r$ produces a homology sphere if and only if $r=1/n$ for some integer $n$. Therefore, the theorem follows from Lemmas~\ref{twist1}, \ref{twist2}, \ref{twist3} and Corollary~\ref{both}.
\end{proof}

\subsection {Proof of Theorem \ref{main_theorem}}

\begin{proof}
A Gabai disk corresponds to a product disk in $\overline{M-S}$. 
As before, one can choose any normal direction of a product disk and
construct a branched surface which corresponds to a taut sutured manifold
decomposition. Then,
either the boundary train track of this branched surface is as shown in
Figure \ref {new.eps} for a certain direction of the Gabai disk, 
or a small neighborhood of the Gabai disk is a 3-ball that contains a
double twist.
If the first case never happens, then Theorem \ref{main_theorem} 
follows from Theorem \ref{Theorem511}.

\begin{figure}[h]
\begin{center}
\resizebox{3in}{!}{\includegraphics{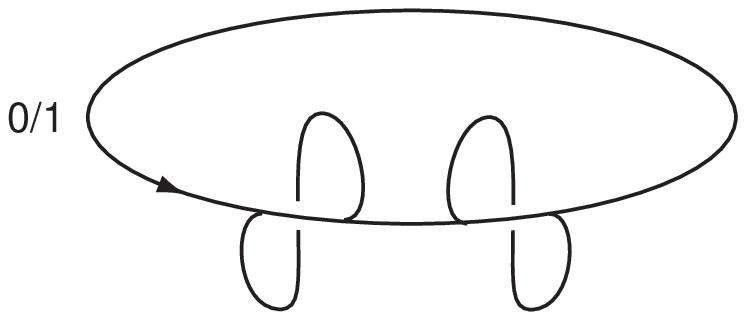}}
\caption{}\label{new.eps}
\end{center}
\end{figure}

Now, suppose the first case happens, i.e. the boundary train track of the
branched surface obtained by adding one Gabai disk is as
shown in Figure \ref{new.eps}. 
Then, similar to the argument above, the branched
surface contains a sink disk iff the two arcs in the boundary of the Gabai
disk are parallel. In this case the knot $K$ is not prime.  
By \cite{DelmanRoberts:PropertyP}, for
a composite knot, any nontrivial Dehn surgery yields a laminar manifold.
If $K$ is prime, then the branched surface contains no sink disk.  By
\cite{TaoLi:TrainTracks}, the train track in Figure \ref{new.eps} 
fully carries any slope in $(-\infty,\infty)$. 
Theorefore, as before, $M(r)$ contains an essential
lamination for any $r\in(-\infty,\infty)$.
\end{proof}

\bibliography{../linklit}

\providecommand{\bysame}{\leavevmode\hbox to3em{\hrulefill}\thinspace}
\providecommand{\MR}{\relax\ifhmode\unskip\space\fi MR }
\providecommand{\MRhref}[2]{%
  \href{http://www.ams.org/mathscinet-getitem?mr=#1}{#2}
}
\providecommand{\href}[2]{#2}
\begin{thebibliography}{CGLS87}

\bibitem[AM90]{AM:Casson}
Selman Akbulut and John~D. McCarthy, \emph{Casson's invariant for oriented
  homology {$3$}-spheres}, Mathematical Notes, vol.~36, Princeton University
  Press, Princeton, NJ, 1990, An exposition. \MR{90k:57017}

\bibitem[BM71]{BingMartin:PropertyP}
R.~H. Bing and J.~M. Martin, \emph{Cubes with knotted holes}, Trans. Amer.
  Math. Soc. \textbf{155} (1971), 217--231. \MR{43 \#4018a}

\bibitem[BM93]{BirmanMenasco:3braids}
Joan~S. Birman and William~W. Menasco, \emph{Studying links via closed braids
  {III}. {C}lassifying links which are closed 3-braids}, Pacific J. Math.
  \textbf{161} (1993), no.~1, 25--113.

\bibitem[Cal00]{Calegari:RFoliations}
Danny Calegari, \emph{The geometry of {${\bf R}$}-covered foliations}, Geom.
  Topol. \textbf{4} (2000), 457--515 (electronic). \MR{2001k:57016}

\bibitem[Cal02]{Calegari:Promoting}
Danny Calegari, \emph{{Promoting Essential Laminations}}, Available as:
  math.GT/0210148, 2002.

\bibitem[CGLS87]{CGLS}
M.~Culler, C.~McA. Gordon, J.~Luecke, and Peter~B. Shalen, \emph{Dehn surgery
  on knots}, Ann. of Math. (2) \textbf{125} (1987), no.~2, 237--300.
  \MR{88a:57026}

\bibitem[DR99]{DelmanRoberts:PropertyP}
Charles Delman and Rachel Roberts, \emph{Alternating knots satisfy {S}trong
  {P}roperty {P}}, Comment. Math. Helv. \textbf{74} (1999), no.~3, 376--397.
  \MR{2001g:57009}

\bibitem[Gab83]{Gabai:FoliationsI}
David Gabai, \emph{Foliations and the topology of $3$-manifolds}, J.
  Differential Geom. \textbf{18} (1983), no.~3, 445--503.

\bibitem[Gab86a]{Gabai:Fibred}
\bysame, \emph{Detecting fibred links in {$S\sp 3$}}, Comment. Math. Helv.
  \textbf{61} (1986), no.~4, 519--555. \MR{88c:57009}

\bibitem[Gab86b]{Gabai:GeneraAlternating}
\bysame, \emph{Genera of the alternating links}, Duke Math. J. \textbf{53}
  (1986), no.~3, 677--681. \MR{87m:57004}

\bibitem[Gab87]{Gabai:FoliationsIII}
\bysame, \emph{Foliations and the topology of {$3$}-manifolds. {III}}, J.
  Differential Geom. \textbf{26} (1987), no.~3, 479--536. \MR{89a:57014b}

\bibitem[GK98]{GK:Laminations}
David Gabai and William~H. Kazez, \emph{Group negative curvature for
  {$3$}-manifolds with genuine laminations}, Geom. Topol. \textbf{2} (1998),
  65--77 (electronic). \MR{99e:57023}

\bibitem[GL89]{GL}
C.~McA. Gordon and J.~Luecke, \emph{Knots are determined by their complements},
  J. Amer. Math. Soc. \textbf{2} (1989), no.~2, 371--415. \MR{90a:57006a}

\bibitem[GO89]{GabaiOertel:EssentialLaminations}
David Gabai and Ulrich Oertel, \emph{Essential laminations in $3$-manifolds},
  Ann. of Math. (2) \textbf{130} (1989), no.~1, 41--73. \MR{90h:57012}

\bibitem[Ker69]{Kervaire:HomologySpheres}
Michel~A. Kervaire, \emph{Smooth homology spheres and their fundamental
  groups}, Trans. Amer. Math. Soc. \textbf{144} (1969), 67--72. \MR{40 \#6562}

\bibitem[KL97]{KL:Genera}
K.~H. Ko and S.~J. Lee, \emph{Genera of some closed 4-braids.}, Topology Appl.
  \textbf{78} (1997), no.~1-2, 61--77.

\bibitem[Li02]{TaoLi:Laminar}
Tao Li, \emph{Laminar branched surfaces in 3-manifolds}, Geom. Topol.
  \textbf{6} (2002), 153--194 (electronic). \MR{1 914 567}

\bibitem[Li03]{TaoLi:TrainTracks}
Tao Li, \emph{Boundary train tracks of laminar branched surfaces}, 2001 Georgia
  International Conference Proceedings, Studies in Advanced Math., vol.~35,
  AMS/IP, 2003, available at: www.math.okstate.edu/$\sim$tli, pp.~269--285.

\bibitem[Lin98]{Lin:FiniteTypeSurvey}
Xiao-Song Lin, \emph{Finite type invariants of integral homology {$3$}-spheres:
  a survey}, Knot theory (Warsaw, 1995), Banach Center Publ., vol.~42, Polish
  Acad. Sci., Warsaw, 1998, pp.~205--220. \MR{99k:57024}

\bibitem[LR02]{LR:CassonHomotopy}
W.~Li and J.~H. Rubinstein, \emph{Casson invariant is a homotopy invariant},
  2002.

\bibitem[MZ00]{MZ:PropertyP}
W.~W. Menasco and X.~Zhang, \emph{Positive knots and knots with braid index
  three have {P}roperty {$P$}}, available as: math.GT/0010154, 2000.

\bibitem[Rob95]{RRoberts:Constructing}
Rachel Roberts, \emph{Constructing taut foliations}, Comment. Math. Helv.
  \textbf{70} (1995), no.~4, 516--545. \MR{96j:57032}

\bibitem[Rob01]{RRoberts:TautFoliation}
\bysame, \emph{Taut foliations in punctured surface bundles. {I}}, Proc. London
  Math. Soc. (3) \textbf{82} (2001), no.~3, 747--768. \MR{2003a:57040}

\bibitem[Sch86]{Scharlemann:PropertyP}
Martin Scharlemann, \emph{A remark on companionship and {P}roperty {P}}, Proc.
  Amer. Math. Soc. \textbf{98} (1986), no.~1, 169--170. \MR{87i:57008}

\bibitem[Sim70]{Simon:PropertyP}
Jonathan Simon, \emph{Some classes of knots with {P}roperty ${P}$}, Topology of
  Manifolds (Proc. Inst., Univ. of Georgia, Athens, Ga., 1969), Markham,
  Chicago, Ill., 1970, pp.~195--199. \MR{43 \#4018b}

\bibitem[ST89]{ST:Linkgenus}
Martin Scharlemann and Abigail Thompson, \emph{Link genus and the {C}onway
  moves}, Comment. Math. Helv. \textbf{64} (1989), no.~4, 527--535.
  \MR{91b:57006}

\bibitem[Sta78]{Stallings:fibred_links}
John~R. Stallings, \emph{Construction of fibred knots and links}, Proc. Symp.
  Pure Math., vol.~32, AMS, 1978, Part 2, pp.~55 -- 59.

\bibitem[Sta96]{Stanford}
Theodore Stanford, \emph{Braid commutators and {V}assiliev invariants}, Pacific
  J. Math. \textbf{174} (1996), no.~1, 269--276. \MR{97i:57008}

\bibitem[Xu92]{Xu}
P.~Xu, \emph{{The Genus of Closed $3$-Braids}}, J. Knot Theory and its Ram.
  \textbf{1} (1992), no.~3, 303 -- 326.

\bibitem[Zha93]{Zhang:PropertyI}
Xingru Zhang, \emph{On property {I} for knots in {$S\sp 3$}}, Trans. Amer.
  Math. Soc. \textbf{339} (1993), no.~2, 643--657. \MR{93m:57011}

\end{thebibliography}
\bibliographystyle {amsalpha}
\end{document}